\documentclass{amsart}

\usepackage{amsmath}
\usepackage{amssymb}
\usepackage{amsthm}
\newtheorem{thm}{Theorem}[section]
\newtheorem{prop}[thm]{Proposition}
\newtheorem{lem}[thm]{Lemma}
\newtheorem{cor}[thm]{Corollary}

\theoremstyle{definition}

\theoremstyle{remark}
\newtheorem{rem}[thm]{Remark}
\newtheorem*{acknowledgement}{Acknowledgements}

\newcommand{\N}{\mathbb{N}}
\newcommand{\R}{\mathbb{R}}
\DeclareMathOperator{\Mag}{Mag}
\DeclareMathOperator{\tr}{tr}
\newcommand{\1}{1\!\!1}
\title[Magnitude and matrix inequalities]
{Magnitude and matrix inequalities}

\author{Kiyonori Gomi}

\address{
Department of Mathematics, 
Institute of Science Tokyo,
2-12-1 Ookayama, Meguro-ku, Tokyo, 152-8551, Japan.}

\email{kgomi@math.titech.ac.jp}

\author{Mark Meckes}

\address{
Department of Mathematics, Applied Mathematics, and Statistics,
Case Western Reserve University
10900 Euclid Ave.
Cleveland, Ohio, 44106,
USA.}

\email{mark.meckes@case.edu}



\begin{document}

\begin{abstract}
  We present several applications of matrix-theoretic inequalities to
  the magnitude of metric spaces.  We first resolve an open problem by
  showing that the magnitude of any finite metric space of negative
  type is less than or equal to its cardinality. This is a direct
  consequence of Styan's matrix inequality involving the Hadamard
  product of matrices. By related methods we also prove a subadditivity
  property for the magnitude function of negative type compact metric
  spaces, and a convexity property for the magnitude for metrics
  interpolating in a natural way between two comparable metrics on a
  given set.
\end{abstract}

\maketitle

\tableofcontents

%%%%%%%%%%%%%%%%%%%%%%%%%%%%%%%%%%%%%%%%%%%%%%%%
%%%%%%%%%%%%%%%%%%%%%%%%%%%%%%%%%%%%%%%%%%%%%%%%

\section{Introduction}
\label{sec:introduction}

The \textit{magnitude} of a metric space is an isometric invariant
introduced by Leinster \cite{Leinster}, which is a formal analogue of the
Euler characteristic of a topological space.  Intuitively, magnitude
is interpreted as the effective number of points in the set, taking
account of the metric. For example, if the distance between two points
is $d$, then the magnitude of this $2$-point metric space is
$2/(1 + e^{-d})$. It is clear that the magnitude goes to $1$ and $2$,
respectively, as the distance $d$ of the two points approaches to $0$
and $\infty$.

We will state the formal definition of the magnitude $\Mag(X)$ of a
metric space $X$, and of some additional terms used in the following
paragraphs, in the next section.  For the moment we will recall that,
for a finite metric space with $n$ points, the definition is based on
properties of an $n \times n$ matrix $Z_X$.  It is thus natural to
suppose that tools from matrix analysis may be used to prove
nontrivial results about magnitude.  However, this avenue has remained
mostly unexplored so far.

In this paper we present three applications to magnitude of
matrix-theoretic inequalities.  The first of these resolves a
fundamental question which has remained open about magnitude since its
introduction, and which crucially supports the interpretation of
magnitude as a metric-sensitive effective number of points.

Since the introduction of magnitude, it has been known that, without
additional assumptions, magnitude can have any real number value, but
with appropriate hypotheses it behaves in a more intuitive way. For
example, if $X$ is a nonempty finite \emph{positively weighted} metric
space, then $1 \le \Mag(X) \le \# X$ if $X$, where $\# X$ means the
cardinality (the number of points) of $X$ \cite[Lemma
2.4.12]{Leinster}. As another assumption, if $X$ is a nonempty finite
\emph{positive definite} metric space, then we have the lower bound
$1 \le \Mag(X)$ \cite[Corollary 2.4.5]{Leinster}. However, positive
definiteness appears to be insufficient to yield an upper bound for
magnitude. An open problem, mentioned by the second author in a
workshop (Magnitude 2023, Osaka), is whether $\Mag(X) \le \#X$ if $X$
is of \textit{negative type}, which is a classical property stronger
than positive definiteness.

The first main result of this paper resolves this problem in the
affirmative.

\begin{thm} \label{thm:main_upper_bound}
If $X$ is a finite metric space of negative type, then $\Mag(X) \le \# X$.
\end{thm}

This upper bound is optimal. Firstly, if $X$ consists of one point,
then $\Mag(X) = \# X = 1$. In addition, let $tX = (X, td)$ denote the
metric space $(X, d)$ rescaled by a real number $t > 0$. The function
$t \mapsto \Mag(tX)$ is known as the \emph{magnitude function} of $X$.
It is known that the magnitude function of any finite metric space is
eventually increasing and tends to $\# X$ as $t \to \infty$
\cite[Proposition 2.2.6]{Leinster}.

After this paper was written, Asao and the first author proved a
stronger version of Theorem \ref{thm:main_upper_bound} by other
methods \cite[Theorem 1.2]{AsGo}.  Namely, they showed that if $X$ is a
finite metric space of negative type and $\#X > 1$, then
$\Mag(X) < \# X$.  Their proof is based on a new geometric
interpretation of magnitude \cite[Theorem 1.1]{AsGo} which was
independently found by Devriendt \cite{Devriendt}.

In this paper, we deduce Theorem \ref{thm:main_upper_bound} as an
immediate consequence of a matrix inequality due to Styan \cite{Styan},
stated as Proposition \ref{prop:Styan} below. Styan's inequality and
the other matrix inequalities applied in this paper are of the form
$A \prec B$, where $A$ and $B$ are symmetric real matrices and $\prec$
denotes the L\"owner order, whose definition will also be recalled in
the next section. The matrices $A$ and $B$ appearing in Styan's
inequality involve the Hadamard (entrywise) product of other
matrices. The Hadamard product arises naturally in the study of
magnitude, and has been explicitly used there before, for example in
\cite[Theorem 3.3]{Meckes}, so our main theorem could have been proved much
earlier. Styan's matrix inequality has been generalized in various
ways, which one may anticipate could imply further interesting
inequalities for magnitude.

As one such application, we prove a subadditivity property for the
magnitude function, although this property can also be given a more
``classical'' magnitude proof.  We also recall that, under the
assumption of positive definiteness (and hence also of negative type),
the notion of magnitude can be extended from finite to compact metric
spaces.

\begin{thm} \label{thm:main_subadditivity}
Let $K$ be a compact metric space of negative type. Then
$$
\Mag((s + t)K) \le \Mag(sK) \Mag(tK)
$$
for every $s, t > 0$. In other words, $\log \Mag(tK)$ is subadditive in $t$.
\end{thm}

Our last main result, Theorem \ref{thm:main_convexity} below, will be
deduced from a different type of matrix inequality, namely, the
operator convexity of the matrix inverse.  Theorem
\ref{thm:main_convexity} involves a different type of interpolation
between two metrics (as opposed to the family of rescalings $td$ for a
given metric $d$), defined in the following lemma.

\begin{lem} \label{lem:metric_interpolation}
Suppose that $X$ is a set equipped with two different metrics $d_0$
and $d_1$ such that $d_0(x,y) \le d_1(x,y)$ for every $x,y \in
X$. Then for each $t \in (0,1)$,
\begin{equation} \label{eq:metric_interpolation}
  d_t(x,y) = - \log \bigl[ (1-t) e^{-d_0(x,y)} + t e^{-d_1(x,y)} \bigr]
\end{equation}
defines a metric on $X$.  Moreover, each $(X, d_t)$ is compact if
$(X,d_1)$ is compact, and is positive definite if both $(X,d_0)$ and
$(X,d_1)$ are positive definite.
\end{lem}

The form of this metric $d_t$ may seem somewhat mysterious.  We will
see in the next section that it actually arises naturally in the
context of the magnitude of metric spaces.

\begin{thm} \label{thm:main_convexity}
  Suppose that $K$ is a set equipped with two metrics $d_0 \le
  d_1$ such that $(K,d_0)$ and $(K,d_1)$ are both compact positive
  definite metric spaces, and let $K_t = (K, d_t)$ for $0\le t \le 1$
  be as defined in \eqref{eq:metric_interpolation}.  Then
  \[
    \Mag(K_t) \le (1-t) \Mag(K_0) + t \Mag(K_1)
  \]
  for every $0 \le t \le 1$.  In other words, $\Mag(K_t)$ is convex in $t$.
\end{thm}

Given the interpretation of magnitude as an effective number of
points, one would intuitively expect that $\Mag(tK)$ would be
increasing in $t$, and that, more generally, $\Mag(K_0) \le \Mag(K_1)$
in the setting of Theorem \ref{thm:main_convexity}.  It is known that
this is false in general; in particular, appropriate rescalings of the
complete bipartite graph $K_{3,2}$ provide a counterexample (see
\cite[Example 2.2.7]{Leinster}). However, this remains an open problem in the
setting of spaces of negative type.  (We note in particular that this
conjecture would also imply Theorem \ref{thm:main_upper_bound}.)
Although Theorem \ref{thm:main_convexity} does not resolve this
problem, the theorem and its proof suggest that the metric
interpolation \eqref{eq:metric_interpolation}, combined with more
powerful matrix inequalities, may be a fruitful approach.

\bigskip

In section \ref{sec:review}, we briefly review  magnitude and the matrix
inequalities used in our proofs. Then in section \ref{sec:proofs}, we
give the proofs of our theorems and discuss some related issues.

\bigskip

\begin{acknowledgement}
  The authors thank Tom Leinster for comments on an earlier version of
  this paper, and the anonymous referee for a thorough reading and
  helpful comments.
\end{acknowledgement}

%%%%%%%%%%%%%%%%%%%%%%%%%%%%%%%%%%%%%%%%%%%%%%%%
%%%%%%%%%%%%%%%%%%%%%%%%%%%%%%%%%%%%%%%%%%%%%%%%
\section{Review of magnitude and matrix inequalities}
\label{sec:review}

%%%%%%%%%%%%%%%%%%%%%%%%%%%%%%%%%%%%%%%%%%%%%%%%
\subsection{The magnitude of finite metric space}

We start with a review of the definition of the magnitude of finite
metric spaces and some relevant notions. All the definitions and the
results here are from Leinster's paper \cite{Leinster} (see also
\cite{LeMe_survey,Meckes}), though our notations and terminology may
differ slightly from the original.

\smallskip

Let $X$ be a finite metric space consisting of $n$ points with $n >
0$.
Without loss of generality, we always realize this as $X = \{ 1,
2, \ldots, n \}$. Suppose that $X$ is equipped with a metric $d : X
\times X \to [0, \infty)$. Then the \textit{zeta matrix} (or the
\textit{similarity matrix}) of the metric space $X = (X, d)$ is
defined as $Z_X = (e^{-d(i, j)}) \in M(n, \R)$. (We write $M(n, \R)$
for the set of $n \times n$ real matrices.)

The form of the zeta matrix of a finite metric space explicates the
relevance to magnitude of both Hadamard products and the metric
interpolation in \eqref{eq:metric_interpolation}.  Observe that the
zeta matrix of the rescaled space $tX$ is $Z_{tX} = (e^{-td(i,j)})$,
from which it follows that $Z_{(s+t)X} = Z_{sX} \circ Z_{tX}$, where
$\circ$ denotes the Hadamard (entrywise) product of matrices.
Moreover, the family of interpolating metrics $\{d_t : 0 < t < 1\}$ in
\eqref{eq:metric_interpolation} is defined precisely so that the
corresponding zeta matrices interpolate linearly in matrix space.
That is,
\begin{equation} \label{eq:Z_t}
  Z_{X_t} = (1-t) Z_{X_0} + t Z_{X_1}.
\end{equation}
It is indeed this relationship that motivates the definition of $d_t$.

A \emph{weighting} on the $n$-point metric space $(X, d)$ is a
vector $w \in \R^n$ such that $Z_X w = \1$, where $\1 \in \R^n$ is the
vector whose entries are $1$. If there exists a weighting
$w$ on $(X, d)$, the magnitude $\Mag(X)$ is defined as the sum of the
entries of $w$, which can be expressed as
\[
\Mag(X) = \langle w, \1 \rangle
\]
using the standard inner product $\langle \cdot ,\cdot \rangle$ on $\R^n$. In
general, weightings need not be unique, but the value of the magnitude
is independent of the choice of a weighting \cite[Lemma 1.1.2]{Leinster}. If
the zeta matrix is invertible, then $X$ has a unique
weighting given by $w = Z_X^{-1} \1$. In this case, the magnitude is
the sum of the entries of $Z_X^{-1}$, which can be expressed as
\begin{equation}
  \label{eq:magnitude_sum}
\Mag(X) = \langle \1, Z_X^{-1} \1 \rangle.
\end{equation}

%For $t > 0$, let $td : X \times X \to [0, \R)$ be the rescaled metric on $X$ 
%defined by $(td)(i, j) = t \cdot d(i, j)$ for $i, j \in X$. 
%The magnitude $\Mag(tX)$ of $tX = (X, td)$, regarded as a function in $t$, 
%is called the magnitude function of $X$.

\bigskip

When a finite metric space $(X, d)$ admits a weighting $w = (w_i)$
whose entries are non-negative, $(X, d)$ is said to be
\textit{positively weighted}. Then, as mentioned in
section \ref{sec:introduction}, the following bounds are known.

\begin{prop}[{\cite[Lemma 2.4.12]{Leinster}}] \label{prop:nnw} Let $X$ be an $n$-point metric
  space with $n > 0$. If $X$ is positively weighted, then
  $1 \le \Mag(X) \le n$.
\end{prop}

\begin{rem} 
  The lower bound given in \cite[Lemma 2.4.12]{Leinster} is $0$ rather than
  $1$, but the improvement above may be known to experts. Considering 
  the first entry of the defining formula $\1 = Z_X w$, one gets
  \[
    1 
    =
    \sum_{j = 1}^n e^{-d(1, j)} w_j 
    \le
    \sum_{j = 1}^n w_j
    = \Mag(X).
  \]
  The lower bound can also be derived from \cite[Lemma 2.2.5]{Leinster}.
\end{rem}

\bigskip

A finite metric space $(X, d)$ is said to be \emph{positive definite}
when $Z_X$ is a positive definite matrix. For a real number $t > 0$,
we write $td : X \times X \to [0, \infty)$ for the rescaled metric on
$X$ given by $(td)(i, j) = t \cdot d(i, j)$, and $tX$ for the rescaled
metric space $(X, td)$. When $Z_{tX}$ is positive definite for all
$t > 0$, the metric space $(X, d)$ is said to be of \emph{negative
  type}. (This classical terminology is due to an equivalent
characterization which we do not explain here.) We then have the
following.

\begin{prop}[{\cite[Corollary 2.4.5]{Leinster}}]
  If $(X, d)$ is a nonempty finite positive definite metric space, then
  $1 \le \Mag(X)$.
\end{prop}

Thus, we also have $1 \le \Mag(X)$ if $(X, d)$ is of negative type.

\bigskip

So far we have been concerned only with finite metric spaces and their
magnitude. The notion of magnitude can be generalized to certain
compact metric spaces \cite{Leinster,LeMe_survey,Meckes}. A compact
metric space $K$ is said to be positive definite (respectively of
negative type) when any finite metric subspace $X \subseteq K$ is
positive definite (respectively of negative type). The magnitude of a
compact positive definite metric space $K = (K, d)$ is defined as
\begin{equation} \label{eq:compact_magnitude}
\Mag(K) = \sup\{ \Mag(X) : X \subseteq K \text{ is finite} \}.
\end{equation}
It follows that $1 \le \Mag(K) \le \infty$ if $K$ is nonempty; it is
known that $\Mag(K) = \infty$ for some compact positive definite
spaces \cite[Theorem 2.1]{LeMe_extremal}.

Note that if $K$ is compact and of negative type, then $\Mag(tK)$ is
defined for every $t > 0$, and $\Mag(tK)$, as a function of $t$, is again
referred to as the \textit{magnitude function} of $K$.

%%%%%%%%%%%%%%%%%%%%%%%%%%%%%%%%%%%%%%%%%%%%%%%%
\subsection{Matrix inequalities}

We first recall the L\"owner order on symmetric  matrices.  If $A, B \in
M(n,\R)$ are symmetric, then we write $A \prec B$ if and only if $B -
A$ is positive semidefinite.  An immediate consequence of this
definition that we will use is that if $A \prec B$, then
\begin{equation} \label{eq:Loewner_sum}
\langle \1, A \1\rangle \le \langle \1, B \1\rangle.
\end{equation}

Styan's matrix inequality that we will use for the proof of Theorem
\ref{thm:main_upper_bound} is the following.

\begin{prop}[{\cite[Corollary 4.2]{Styan}}] \label{prop:Styan}
Let $R \in M(n, \R)$ be a positive definite symmetric matrix whose diagonal entries are $1$. Then
$$
2(R \circ R)^{-1} \prec R^{-1} \circ R + I,
$$
where $I \in M(n, \R)$ denotes the identity matrix.
\end{prop}

In Styan's paper \cite{Styan}, Proposition \ref{prop:Styan} is proved
from the viewpoint of multivariate statistical
analysis. Matrix-theoretic proofs (of generalizations) can be found in
\cite{Ando,Liu,Visick,WaZh,YaFeCh,Zhang} for instance.

It should be noted that $R \circ R$ is positive definite, thanks to
the classical Schur product theorem, which we state here for
reference.

\begin{prop}[see, e.g., {\cite[Theorem 3.1]{Styan}}] \label{prop:Schur}
If $A, B \in M(n,\R)$ are both positive definite, then $A \circ B$ is
also positive definite.
\end{prop}

\medskip

As a direct consequence of Proposition \ref{prop:Styan}, we obtain the
following.

\begin{lem} \label{lem:inequality}
For any positive definite symmetric matrix $R \in M(n, \R)$ whose
diagonal entries are $1$, we have
$$
\langle \1, (R \circ R)^{-1} \1 \rangle \le n.
$$
\end{lem}

\begin{proof}
Proposition \ref{prop:Styan} and \eqref{eq:Loewner_sum} imply that
\begin{align*}
2\langle \1, (R \circ R)^{-1} \1 \rangle
&\le \langle \1, (R^{-1} \circ R) \1 \rangle + n
=
\tr(R^{-1} R) + n
= 2n,
\end{align*}
where we have also used the elementary fact that
$\langle \1, (A \circ B^t) \1 \rangle = \tr (AB)$ (see \cite[formula (2.8)]{Styan}).
\end{proof}

\medskip

In one of our proofs of Theorem \ref{thm:main_subadditivity}, we will
apply the following Styan-type inequality due to Wang and Zhang.

\begin{prop}[{\cite[Theorem 1]{WaZh}}] \label{prop:Wang_Zhang}
For any positive definite symmetric matrices $A, B \in M(n, \R)$ and
any matrices $C, D \in M(m, n, \R)$, we have
$$
(C \circ D)(A \circ B)^{-1}(C \circ D)^t \prec
(CA^{-1}C^t) \circ (DB^{-1}D^t).
$$
\end{prop}

\medskip

Our proof of Theorem \ref{thm:main_convexity} uses the following
theorem, which says that the matrix inverse is an operator convex
mapping on the space of positive definite matrices; see for example
\cite[Corollary 4.1]{Ando}.

\begin{prop} \label{prop:inverse_convex}
  Let $A,B \in M(n,\R)$ be
  positive definite.  Then for every $0 \le t \le 1$, we have
  \[
    \bigl[(1-t) A + t B]^{-1} \prec (1-t) A^{-1} + t B^{-1}.
  \]
\end{prop}

%%%%%%%%%%%%%%%%%%%%%%%%%%%%%%%%%%%%%%%%%%%%%%%%
%%%%%%%%%%%%%%%%%%%%%%%%%%%%%%%%%%%%%%%%%%%%%%%%
\section{Proofs of the main theorems}
\label{sec:proofs}

%%%%%%%%%%%%%%%%%%%%%%%%%%%%%%%%%%%%%%%%%%%%%%%%
\subsection{The upper bound}

We now prove our first main theorem, Theorem
\ref{thm:main_upper_bound}, in an equivalent form.

\begin{thm} \label{thm:main_upper_bound_repeated} Let $X$ be an
  $n$-point metric space. If $X$ is of negative type, then we
  have $\Mag(tX) \le n$ for any $t > 0$.
\end{thm}

\begin{proof}
  By the assumption that $X$ is of negative type, the symmetric matrix
  $Z_{\frac{t}{2}X} = (e^{-\frac{t}{2}d(i, j)})$ is positive definite
  for all $t > 0$. Noting that
  $Z_{\frac{t}{2}X} \circ Z_{\frac{t}{2}X} = Z_{tX}$ and applying
  Lemma \ref{lem:inequality} to $R = Z_{\frac{t}{2}X}$, we have by
  \eqref{eq:magnitude_sum} that
$$
\Mag(tX) = \langle \1, Z_{tX}^{-1} \1 \rangle
= \langle \1, (Z_{\frac{t}{2}X} \circ Z_{\frac{t}{2}X})^{-1} \1 \rangle
\le n,
$$
as claimed.
\end{proof}

In view of the argument in the proof of Theorem
\ref{thm:main_upper_bound_repeated}, the assumption that $(X, d)$ is
of negative type can be weakened as follows.

\begin{thm} \label{thm:weak_assumption} Let $X = (X, d)$ be an
  $n$-point metric space. If $\frac{1}{2}X = (X, \frac{1}{2}d)$ is
  positive definite, then $\Mag(X) \le n$.
\end{thm}

If a finite metric space $X$ is positively weighted and $\# X \ge 2$,
then the proof of Proposition \ref{prop:nnw} in \cite{Leinster} also shows
that $\Mag(X) < \# X$. (Note that we demand $d(i, j) < \infty$ for $d$
to be a metric.) This applies in particular to spaces with 
$\# X \le 3$ \cite[Proposition 2.4.15]{Leinster}. The proof of
Theorem \ref{thm:main_upper_bound_repeated}, shows that
$\Mag(X) = \# X$ if and only if the equality in Styan's matrix
inequality is attained.  Equality conditions in Styan's inequality are
given in \cite{YaFeCh} for example, but these are not easily
interpretable in terms of the geometry of $X$.

The results here and the discussion above leave open the question of
whether $\Mag(X) < \#X$ for all finite spaces of negative type with
$\# X \ge 4$. As mentioned in the introduction, Asao and the first
author have since proved this using a new geometric interpretation of
magnitude \cite[Theorem 1.1]{AsGo}.

%%%%%%%%%%%%%%%%%%%%%%%%%%%%%%%%%%%%%%%%%%%%%%%%
%%%%%%%%%%%%%%%%%%%%%%%%%%%%%%%%%%%%%%%%%%%%%%%%
\subsection{A subadditivity property}
\label{subsec:subadditivity}

We next prove Theorem \ref{thm:main_subadditivity} in two ways.  The
first is based directly on ``classical'' magnitude techniques.

\begin{proof}[First proof of Theorem \ref{thm:main_subadditivity}]
  We recall the tensor product (or $\ell^1$ sum) of two metric spaces
  $A$ and $B$ (see \cite[Example 1.4.2]{Leinster}): $A \otimes B$ is the set
  $A \times B$ equipped with the metric
  \[
    d((a,b), (a', b')) = d(a,a') + d(b,b').
  \]
  Proposition 3.1.4 of \cite{Leinster} implies that if $A$ and $B$ are
  positive definite then $A \otimes B$ is as well, and that
  $\Mag(A\otimes B) = \Mag(A) \Mag(B)$; hence in particular we have
  \begin{equation} \label{eq:tensor-product}
    \Mag((sK) \otimes (tK)) = \Mag(sK) \Mag(tK).
  \end{equation}

  The definition \eqref{eq:compact_magnitude} implies that magnitude
  is weakly increasing with respect to set inclusion for compact
  positive definite metric spaces.  The diagonal map $x \mapsto (x,x)$
  defines an isometric embedding of $(s+t) K$ into
  $(sK) \otimes (tK)$, which implies that
  \begin{equation}\label{eq:log-subadditivity}
    \Mag((s+t)K)
    \le \Mag((sK) \otimes (tK)).
  \end{equation}
  Combining \eqref{eq:tensor-product} and \eqref{eq:log-subadditivity}
  proves the claim.
\end{proof}

Theorem \ref{thm:main_subadditivity} can also be given a quicker
(albeit less elementary) proof similar in spirit to the proof of
Theorem \ref{thm:main_upper_bound} via Styan's inequality.

\begin{proof}[Second proof of Theorem \ref{thm:main_subadditivity}]
%Theorem 1 of \cite{WZ} states that, if $A,B \in M(n,\R)$
%are positive definite and $C,D \in M(m,n,\R)$, then
%\begin{equation}
%  \label{eq:Wang-Zhang}
%  (C A^{-1} C^t) \circ (D B^{-1} D^t) - (C \circ D) (A\circ B)^{-1}
%  (C\circ D)^t
%\end{equation}
%is positive semidefinite.
Suppose first that $K$ is finite with $n$ points.  Applying
%\eqref{eq:Wang-Zhang} 
Proposition \ref{prop:Wang_Zhang}
with $A = Z_{sK}$, $B = Z_{tK}$, and
$C = D = \1^t \in M(1,n,\R)$ immediately implies the theorem in this
case.  In the general case, for each finite subset $X \subseteq K$, we
then have
\[
  \Mag((s+t)X) \le \Mag(sX) \Mag(tX) \le \Mag(sK) \Mag(tK)
\]
by \eqref{eq:compact_magnitude}. The general case now follows since
\[
  \Mag((s+t)K) = \sup\{ \Mag((s+t)X) : X \subseteq K \text{ is finite}\}.
\qedhere
\]
\end{proof}

As with Theorems \ref{thm:main_upper_bound_repeated} and Theorem
\ref{thm:weak_assumption}, both of the proofs of Theorem
\ref{thm:main_subadditivity} yield a more general result. Note that
the positive definiteness of $(K, d_1 + d_2)$ follows from the Schur
product theorem (Proposition \ref{prop:Schur}).

\begin{thm}
  \label{thm:subadditivity}
  Let $K$ be a set equipped with two metrics $d_1$ and $d_2$ such that
  $(K, d_1)$ and $(K, d_2)$ are both compact and positive
  definite. Define a new metric $d_1 + d_2$ on $K$ by
  $(d_1 + d_2)(x, y) = d_1(x, y) + d_2(x, y)$. Then $(K, d_1+d_2)$ is
  positive definite and
  \[
    \Mag(K, d_1 + d_2) \le \Mag(K, d_1) \Mag(K, d_2).
  \]
  
  In particular, if $(K,d)$ is a compact metric space such that
  $s, t > 0$ and $sK$ and $tK$ are positive definite, then $(s+t)K$ is
  positive definite and $\Mag((s+t)K) \le \Mag(sK) \Mag(tK)$.
\end{thm}

An easy induction argument proves the following previously known
consequence of Theorem \ref{thm:subadditivity}.

\begin{cor}[{\cite[Proposition 3.3]{LeMe_extremal}}]
  \label{cor:subexponential}
  Let $X$ be a compact positive definite metric space.  Then $\Mag(nX)
  \le \Mag(X)^n$ for every $n \in \N$.
\end{cor}

Corollary \ref{cor:subexponential} is a step in the direction of a
conjecture made by Mokshay Madiman (personal communication).  Suppose
that $E$ is a positive definite finite-dimensional normed space.
(Equivalently, $E$ is isometrically isomorphic to a finite-dimensional
subspace of $L^1$; see \cite[Corollary 3.5]{Meckes}.)  Madiman conjectured the following
Brunn--Minkowski inequality for magnitude.  If $X, Y \subseteq E$ are
compact (possibly satisfying additional hypotheses) and $0 \le \lambda
\le 1$, then
\begin{equation}
  \label{eq:Brunn-Minkowski}
  \Mag(\lambda X + (1-\lambda) Y) \ge \Mag(X)^\lambda
  \Mag(Y)^{1-\lambda}.
\end{equation}
Here $+$ refers to Minkowski addition of sets in the vector space $E$.

This conjecture is known to be true in this generality in one
dimension (see \cite[Theorem 1.9]{ALMM}).  Thanks to \cite[Theorem
4.1.4]{LeMe_survey}, \eqref{eq:Brunn-Minkowski} would imply the classical
Brunn--Minkowski inequality, where magnitude is replaced by Lebesgue
measure.  Unlike the classical case, however,
\eqref{eq:Brunn-Minkowski} would be a nontrivial statement even when
$Y = \{0\}$, in which case it reduces to the statement
\begin{equation}
  \label{eq:BM_one_set}
  \Mag(\lambda X) \ge \Mag(X)^\lambda,
\end{equation}
which makes sense for arbitrary spaces of negative type.  Substituting
$\frac{1}{n} X$ in place of $X$, Corollary \ref{cor:subexponential}
implies that \eqref{eq:BM_one_set} holds for $\lambda = \frac{1}{n}$.

%%%%%%%%%%%%%%%%%%%%%%%%%%%%%%%%%%%%%%%%%%%%%%%%
%%%%%%%%%%%%%%%%%%%%%%%%%%%%%%%%%%%%%%%%%%%%%%%%
\subsection{A convexity property}
\label{subsec:convexity}

Finally, we prove Lemma \ref{lem:metric_interpolation} and our last
main result, Theorem \ref{thm:main_convexity}.

\begin{proof}[Proof of Lemma \ref{lem:metric_interpolation}]
  The fact that $d_t(x,y) \ge 0$ with equality if and only if $x=y$
  follows immediately from the fact that $d_0$ and $d_1$ are metrics
  and the fact that the exponential function and logarithm are
  strictly increasing.  The
  triangle inequality for $d_t$ is equivalent to the statement that
  \begin{multline} \label{eq:multiplicative_triangle_inequality}
    \bigl[(1-t) e^{-d_0(x,y)} + t e^{-d_1(x,y)}\bigr]
    \bigl[(1-t) e^{-d_0(y,z)} + t e^{-d_1(y,z)}\bigr] \\
    \le
    (1-t) e^{-d_0(x,z)} + t e^{-d_1(x,z)}
  \end{multline}
  for every $x,y,z \in X$.  Given $x,y,z \in X$, denote
  $\alpha = e^{-d_0(x,y)}$, $\beta = e^{-d_1(x,y)}$,
  $\gamma = e^{-d_0(y,z)}$, and $\delta = e^{-d_1(y,z)}$.  The
  triangle inequalities for $d_0$ and $d_1$ imply that 
  \eqref{eq:multiplicative_triangle_inequality} would follow if
  \[
    \bigl[(1-t) \alpha + t \beta\bigr] \bigl[(1-t) \gamma + t\delta
    \bigr]
    \le (1-t) \alpha \gamma + t \beta \delta.
  \]
  By hypothesis we have that $\beta \le \alpha$ and $\delta \le
  \gamma$, and therefore
  \[
    (1-t) \alpha \gamma
    + t \beta \delta - \bigl[(1-t) \alpha
    + t \beta\bigr] \bigl[(1-t) \gamma + t\delta \bigr]
    = t (1-t) (\alpha - \beta)(\gamma - \delta) \ge 0,
  \]
  which completes the proof that $d_t$ is a metric.

  Since $d_t \le d_1$ for all $0 \le t \le 1$, the identity map from
  $(X,d_1)$ to $(X,d_t)$ is continuous, and hence $(X,d_t)$ is compact
  if $(X,d_1)$ is.

  Finally, to prove the statement about
  positive definiteness it suffices to assume that $X$ is finite with
  $n$ points.  The claim then follows from the fact that $Z_{X_t}$ is
  positive definite by \eqref{eq:Z_t} and the convexity of the set of
  positive definite $n\times n$ matrices.
\end{proof}

\begin{proof}[Proof of Theorem \ref{thm:main_convexity}]
  As in the second proof of Theorem \ref{thm:main_subadditivity}, it
  suffices to assume that $X$ is finite with $n$ points.  Proposition
  \ref{prop:inverse_convex} and \eqref{eq:Z_t} imply that
  \[
    Z_{X_t}^{-1} \prec (1-t) Z_{X_0}^{-1} + t Z_{X_1}^{-1},
  \]
  which by \eqref{eq:magnitude_sum} and \eqref{eq:Loewner_sum} proves
  the theorem.
\end{proof}

%%%%%%%%%%%%%%%%%%%%%%%%%%%%%%%%%%%%%%%%%%%%%%%%
%%%%%%%%%%%%%%%%%%%%%%%%%%%%%%%%%%%%%%%%%%%%%%%%

\bibliographystyle{abbrv}
\bibliography{magnitude_matrix_inequalities}

\end{document}